# An Augmented Smoothing Method of $L_1$-norm Minimization and Its Implementation by Neural Network Model

Yunchol Jong [a]

[a] Center of Natural Sciences, University of Science, Pyongyang, DPR Korea
E-mail: yuncholjong@yahoo.com

**Abstract.** In this paper we propose an augmented smoothing function for nonlinear $L_1$–norm minimization problem and consider a global stability of a gradient-based neural network model to minimize the smoothing function. The numerical simulations show that our smoothing neural network finds successfully the global solution of the $L_1$-norm minimization problems considered in the simulation.

## 1. Introduction

The nonlinear $L_p$-norm estimation problem is to find $x$ which minimizes

$$\min f(x) = \sum_{i=1}^{m} |r_i(x)|^p, \ x \in R^n,$$

where $1 \leq p \leq \infty$ and $r_i : R^n \to R, i = 1,\ldots, m$ are smooth nonlinear functions of $x$. In data-fitting problems, $r_i$ are residual functions.

The discrete $L_p$-approximation problem is a basic problem in approximation theory and optimization. This problem is of basic importance in approximation theory with applications in data processing, engineering design, operations research, and other fields of applied sciences. In particular, $L_1-$ approximation problem has been widely applied to the approximation, statistical processing and circuit optimization because the robustness of the $L_p$ estimates becomes stronger when p tends to 1. The $L_p$- problem could be solved either by transforming it into a nonlinear programming problem, or as a general unconstrained minimization problem, but none of them could be better than algorithms which take the special structure of the $L_p$- problem into account. All the algorithms were inspired by nonlinear least-squares methods. The $L_1$-results could be got through the algorithms of $L_p$-problem by letting $p \to 1$.

Most algorithms that have been proposed for solving this problem are essentially variants of Newton's method. However, since the objective function is not continuously differentiable for $p=1$, the Newton methods need essential revisions to handle this case. [10], [11] and [12] considered discrete linear L₁-approximation whose objective function is $\sum_{i=1}^{m} |y_i - a_i^T x|$. In [5] was proposed a method of analytic centers for minimization of $\sum_{i=1}^{m} |b_i + a_i^T x|^p$, p ≥1 and in [4]

was proposed an algorithm for nonlinear $L_p$ - norm minimization.

In the last two decades considerable attention has been focused on optimization neural networks. Such systems are considered as potentially efficient hardware solutions for large-scale or hard optimization problems. Optimization neural networks could work very fast as a parallel computational structure in truly distributed implementation. The gradient dynamical system can be solved using standard ordinary differential equation (ODE) software and this could be an advantage, when the number of unknowns is large. Furthermore, implementations of these standard ODE methods on parallel computers could be used, in order to make it possible to deal with large datasets. In addition, gradient dynamical systems can be implemented as an analog circuit using only resistors, amplifiers and switches, which is appropriate for real time processing using VLSI technology [3]. In [3] was proposed a dynamical gradient system with discontinuous right-side to find a solution set of linear $L_1$ - norm problem. They proved finite time convergence to the solution set by means of a diagonal type Lyapunov function for the dynamical system. In [1] and [2] was proposed a discrete-time neural network for solving constrained linear $L_1$ estimation problems fast.

To author's knowledge, there is no any gradient-based neural network for nonlinear $L_1$ - norm problem. In [6] was presented a neural network using eigenvalues and eigenvectors of Hessian matrix of the objective function to find a local minimum of the problem. Their method should continue to find equilibrium points of the neural network to obtain an equilibrium point satisfying second-order optimality condition, and need a lot of computation.

In this paper we propose an augmented smoothing function for nonlinear $L_1$ - problem and consider a global stability of a gradient-based neural network model to minimize the smoothing function. The numerical simulations show that our smoothing neural network finds successfully the global solution of the $L_1$ - problem.

## 2. An augmented smoothing for nonlinear $L_1$ - problem

The nonlinear $L_1$ -problem is defined as follows.

$$\min f(x) = \sum_{i=1}^{m} |r_i(x)|, x \in R^n, \qquad (1)$$

where $f(x)$ can be considered to be $L_1$ -norm of residual vector $r(x) = (r_1(x),...,r_m(x))^T$ and $r_i : R^n \to R, \quad i=1,...,m$ are twice continuously differentiable.

Let $|r_i(x)|$ be approximated by the following smooth function.

$$s_i(x,\mu) = \mu^2 \ln\left[\exp\left(\frac{r_i(x)}{\mu^2}\right) + \exp\left(-\frac{r_i(x)}{\mu^2}\right)\right], \mu \neq 0, \qquad (2)$$

where $\mu$ is a smoothing variable. We consider $\mu$ a variable as well as $x$, while $\mu$ is considered to be a parameter in usual smoothing methods.

**Lemma 1**. We have the following relations between $s_i(x,\mu)$ and $|r_i(x)|$.

(i) $|r_i(x)| < s_i(x,\mu) < |r_i(x)| + \mu^2 \ln 2, \ \mu \neq 0$ \hspace{2em} (3)

$$s_i(x,0) = |r_i(x)| \tag{4}$$

(ii) When $\mu \neq 0$, $s_i(x, \mu)$ is continuously differentiable and

$$\nabla_x s_i(x,\mu) = \alpha_i(x,\mu)\nabla r_i(x), \tag{5}$$

$$\alpha_i(x,\mu) = \frac{\exp\left(\frac{2r_i(x)}{\mu^2}\right) - 1}{\exp\left(\frac{2r_i(x)}{\mu^2}\right) + 1}, \quad |\alpha_i(x,\mu)| < 1 \tag{6}$$

$$sign(\alpha_i(x,\mu)) = sign(r_i(x))$$

$$\alpha_i(x,\mu) \to 1(\mu \to 0, r_i(x) > 0), \ \alpha_i(x,\mu) \to -1(\mu \to 0, r_i(x) < 0), \tag{7}$$

$$\frac{\partial s_i(x,\mu)}{\partial \mu} = \frac{2}{\mu}(s_i(x,\mu) - \alpha_i(x,\mu)r_i(x)) \tag{8}$$

(iii) When $\mu \neq 0$, $s_i(x,\mu)$ is twice continuously differentiable and

$$\nabla_x^2 s_i(x,\mu) = \frac{2\lambda_i(x,\mu)}{\mu^2}(1-\alpha_i(x,\mu))\nabla r_i(x)\nabla r_i(x)^T + \alpha_i(x,\mu)\nabla^2 r_i(x)$$

$$\nabla_{x\mu}^2 s_i(x,\mu) = -\frac{4\lambda_i(x,\mu)}{\mu^3}r_i(x)(1-\alpha_i(x,\mu))\nabla r_i(x) \tag{9}$$

$$\frac{\partial^2 s_i(x,\mu)}{\partial \mu^2} = \frac{2}{\mu^2}(s_i(x,\mu) - \alpha_i(x,\mu)r_i(x)) + 2\left(\frac{2r_i(x)}{\mu^2}\right)^2 (1-\alpha_i(x,\mu))\lambda_i(x,\mu),$$

where

$$\lambda_i(x,\mu) = \frac{\exp\left(\frac{2r_i(x)}{\mu^2}\right)}{\exp\left(\frac{2r_i(x)}{\mu^2}\right) + 1} \tag{10}$$

(Proof) The proposition (i) is the immediate consequence of the Lemma 1 of [6]. The relations of proposition (ii) and (iii) are obtained by direct calculation. □

Now, let $z = (x, \mu)$, $E(z) = \sum_{i=1}^{m} s_i(z)$ and $E_1(z) = E(z) + \frac{\mu^2 \|x\|^2}{2}$. And let

$$\partial^0 E_1(z) = \begin{cases} \{\nabla E_1(z)\}, & \mu \neq 0 \\ \{\lim_{\mu \to 0} \nabla E_1(z)\}, & \mu = 0 \end{cases},$$

where $\nabla E_1(z) = \begin{pmatrix} \nabla_x E_1(z) \\ \frac{\partial E_1(z)}{\partial \mu} \end{pmatrix}$.

**Theorem 1.** If $0 \in \partial^0 E_1(z)$, then $\mu = 0$ and $x$ is a stationary point of $f(x)$. Moreover, if $r_i(x), i = 1,...,m$ are affine functions, then $x$ is a solution of the problem (1).

(Proof) If $\mu \neq 0$, we have

$$\alpha_i(x,\mu)r_i(x) \leq |r_i(x)| < s_i(x,\mu), \quad i = 1,\ldots,m \quad (11)$$

by (6) and (3). In view of (5), we have

$$\frac{\partial E_1(z)}{\partial \mu} = \frac{2}{\mu}\sum_{i=1}^{m}(s_i(x,\mu) - \alpha_i(x,\mu)r_i(x)) + \mu\|x\|^2 \quad (12)$$

and $s_i(x,\mu) - \alpha_i(x,\mu)r_i(x) > 0, i = 1,\ldots,m$ by (11). Therefore, $\dfrac{\partial E_1(z)}{\partial \mu} \neq 0$ for $\mu \neq 0$, which implies that $\nabla E_1(z) \neq 0$. Hence $0 \in \partial^0 E_1(z)$ implies that $\mu = 0$, and $\alpha_i(x,0) = sign(r_i(x))$ for $i \in \{j \mid r_j(x) \neq 0\}$ and $\alpha_i(x,0) \in [-1,1]$ for $i \in \{j \mid r_j(x) = 0\}$ by (5)-(7). Thus $0 \in \partial^0 E_1(z)$ implies that

$$0 \in \left\{\lim_{\mu \to 0} \nabla E_1(z)\right\}.$$

Since

$$\lim_{\mu \to 0} \nabla E_1(z) = \sum_{i=1}^{m} \alpha_i(x,0)\nabla r_i(x) \in \partial f(x), \quad (13)$$

We have $0 \in \partial f(x)$, i.e., $x$ is a stationary point of $f(x)$, where $\partial f(x)$ is the subdifferential of $f(x)$ at $x$ defined by

$$\partial f(x) = \left\{a \mid a = \sum_{i=1}^{m} \delta_i \nabla r_i(x)\right\},$$

where $\delta_i = sign(r_i(x))$ if $r_i(x) \neq 0$ and $\delta_i \in [-1, 1]$ if $r_i(x) = 0$.

Besides, if $r_i(x), i = 1,\ldots,m$ are affine functions, then $f(x)$ is convex function and the stationary point $x$ is a solution of the problem (1). □

From the above theorem, we see that a solution of the problem (1) can be obtained by finding a stationary point of the function $E_1(z)$.

## 3. A neural network model for the augmented smooting

Let us consider the following neural network model for finding a stationary point of $E_1(z)$.

$$\frac{dz}{dt} = -M\nabla E_1(z), \quad (14)$$

where $M$ is a positive diagonal matrix for scaling. In what follows, let $M = \tau I$, where $\tau$ is a positive constant and $I$ is an identity matrix.

**Theorem 2.** Suppose that $\nabla E_1(z)$ satisfies the Lipschitz condition. Then the trajectory $z(t)$

of the dynamical system (14) satisfying the initial condition $z^0 = z(0)$ for arbitrary $z^0$ converges to an equilibrium point $z^* = (x^*, \mu^*)$ as $t \to \infty$, and $\mu^* = 0$ and $x^*$ is a stationary point of the problem (1).

(Proof) There exists a unique solution $z(t)$ of the initial problem (14) satisfying $z^0 = z(0)$ for arbitrary $z^0$ in $[0, \infty)$ because the right-side of (14) satisfies the Lipschitz condition. In view of (3), we have that

$$f(x) \leq E(z) \leq f(x) + \mu^2 m \ln 2, \ f(x) \geq 0$$

for arbitrary $x$ and $z = (x, \mu)$ and $E_1(z) \geq 0$ for every $z$, i.e., $E_1(z)$ is bounded below. Thus, the trajectory $z(t)$ of the system (14) converges to an equilibrium point $z^* = (x^*, \mu^*)$ as $t \to \infty$ by Theorem 4 of [4], and $\mu^* = 0$ and $x^*$ is a stationary point of the problem (1) by Theorem 1. □

**Corollary 1**. Suppose that $\nabla E_1(z)$ satisfies the Lipschitz condition and $E_1(z) \to +\infty$ as $\|z\| \to \infty$. If the equilibrium point of the system (13) is unique, then the equilibrium point is globally and asymptotically stable.

(Proof) By virtue of the assumption, there exists a unique solution $z(t)$ of (14) satisfying the initial condition for any $z^0 = z(0)$ in time interval $[0, \infty)$.
Then

$$\frac{dE_1(z)}{dt} = \nabla E_1(z)^T \frac{dz}{dt} = -\tau \nabla E_1(z)^T \nabla E_1(z) = -\tau \|\nabla E_1(z)\|^2 \leq 0 \tag{15}$$

and $W(z) = \tau \|\nabla E_1(z)\|^2 > 0$ for any $z \neq z^*$ if the equilibrium point $z^*$ is unique. Since $E(z)$ is bounded below, we have global and asymptotic stability of the equilibrium point $z^*$ by Theorem 2 of [5]. □

## 4. Numerical simulation

Our algorithms have been implemented by 1.5GHZ 20GB 256MB PC NEC in the Matlab7.8 environment using ode23. In numerical experiment, we used the following smoothing function instead of (2):

$$s_i(x, \mu) = \theta \mu^2 \ln \left[ \exp\left(\frac{r_i(x)}{\theta \mu^2}\right) + \exp\left(-\frac{r_i(x)}{\theta \mu^2}\right) \right]$$

where $\theta$ is a positive constant.
**Problem 1.** We carried out numerical simulation for the following problem considered in [4].

$$\min f(x) = \left(x_1^3 - 3x_1\right)^2 + x_2^2 \tag{16}$$

This problem has five stationary points $(0,0), (\pm\sqrt{3}, 0)$ and $(\pm 1, 0)$ in which three points $(0,0)$

and $(\pm\sqrt{3},0)$ are optimal solutions.
The simulation result of [4] is the following:

Table 1.

| Starting point | $t_f$ | $\nabla E(x)$ | $E(x)$ | Limit point |
|---|---|---|---|---|
| (1,1) | 0.0107 | 5.87e-07 | 2.43e-15 | (1.73205,0) |
| (-1,-1) | 0.0079 | 5.51e-07 | 9.48e-15 | (-3.04e-08,3.43e-08) |
| (2,2) | 0.0076 | 6.11e-07 | 9.35e-14 | (1.73205,0) |
| (-2,-2) | 0.0076 | 6.11e-07 | 9.35e-14 | (-1.73205,0) |

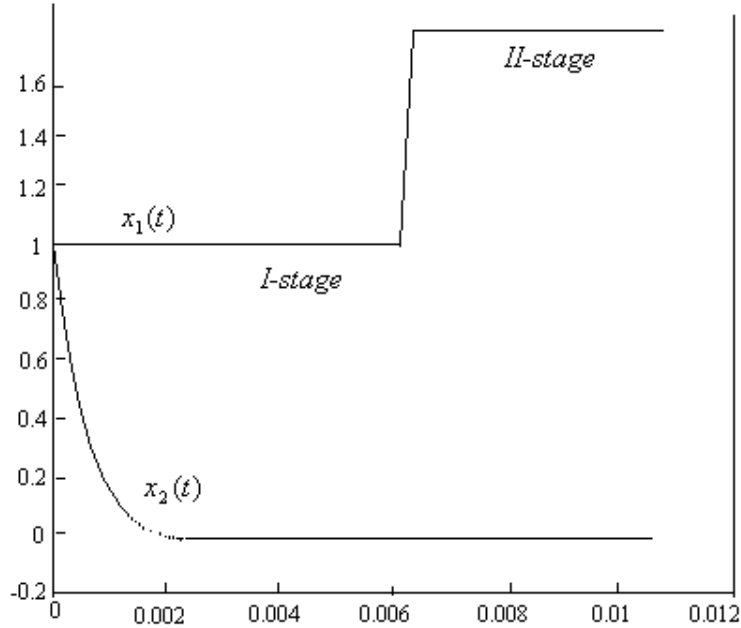

Fig.1. Behavior of the trajectory of [4] starting at $x^0 = (1,1)$

The problem is nonconvex and one solution was obtained by applying neural network approach twice in [4](see Fig.1). When starting at points (1,1) and (-1,-1), trajectories of their neural network converged to point (1,0) and (-1,0), respectively, which are stationary points but not solutions. Thus, they carried out simulation again by taking these points as starting points and obtained the result shown in the above table, where $t_f$ is end of the integration interval.

The problem (16) is equivalent to the following problem.

$$\min\ f(x) = |x_1^3 - 3x_1| + |x_2|$$

We applied our neural network approach of augmented smoothing (NNAS) and neural network approach (NN) proposed in [1] to the above problem, respectively.
When letting

$$M = \begin{pmatrix} 10000 & 0 & 0 \\ 0 & 10 & 0 \\ 0 & 0 & 100 \end{pmatrix},\ \mu^0 = 40,\ \theta = 0.02,$$

simulation result is as following (Figure 2, 3 and Table 2):

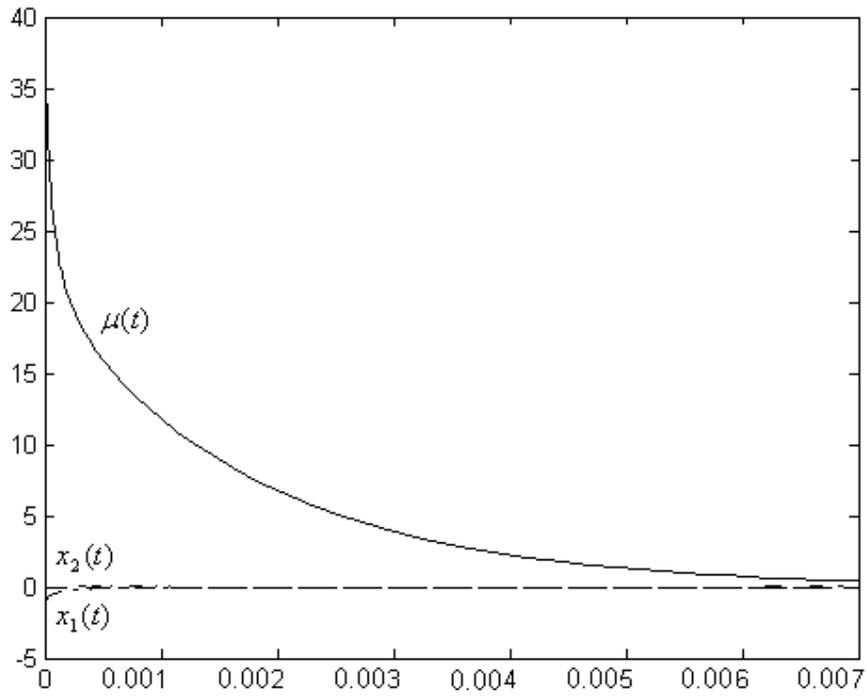

Fig.2. Behavior of our trajectory starting at $x^0 = (40, -1, -1)$

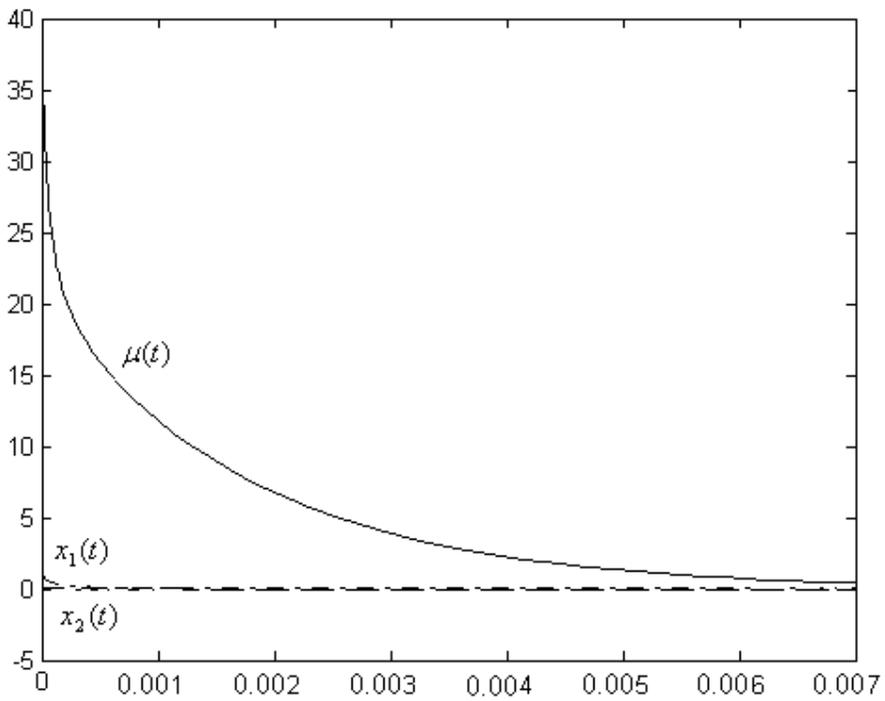

Fig.3. Behavior of our trajectory starting at $x^0 = (40, 1, 1)$

Table 2.

| Starting point | NNAS in [0,0.007] | | | NN in [0,0.01] | | |
|---|---|---|---|---|---|---|
| | time | Object value | solution | time | Object value | solution |
| (1, 1) | 0.25 | 8.6711e-07 | 1.0e-06 * (0.2890, -0.0000) | 0.219 | 2.0000 | (1.0000, 0.0000) |
| (-1, -1) | 0.266 | 8.6711e-07 | 1.0e-06 * (-0.2890, 0.0000) | 0.094 | 2.0000 | -(1.0000, 0.0000) |
| random (50) | 0.203 | 6.5288e-08 | 1.0e-07 * (-0.2176, 0.0000) | 76.016 (7) | 1.4765e-06 | 1.0e-06 * (-0.0204, 0.7876) |
| | | | | 95.45 (32) | 1.8042 | (0.7323, 0.0000) |
| | | | | 150.64 (11) | 1.2322 | (-0.4389, 0.0000) |

As shown in the table and figures, our approach NNAS always converged to the solution of the given problem with any starting point of [-1,1]×[-1,1] and is better in speed and accuracy than NN of [1]. In the table, the number in round brackets of the time column is the frequency of success among runs.

**Problem 2**(Rastrigin function). The following problem is equivalent to Rastrigin's multiextremal problem.

$$\min f(x) = \left|x_1^2 + x_2^2\right| + \left|20 - 10 \cdot [\cos(2\pi x_1) + \cos(2\pi x_2)]\right|,$$
$$\text{s.t.} \quad -1 \leq x_i \leq 1, \ i = 1, 2$$

The solution of the problem is (0, 0).

With the same $M, \mu^0$ and $\theta$ as in problem 1, simulation results are as following.

Table 3.

| starting point | NNAS, $t_f = 0.007$ | | | NN, $t_f = 0.01$ | | |
|---|---|---|---|---|---|---|
| | time | $f(x)$ | $x$ | time | $f(x)$ | $x$ |
| -(1,1) | 0.232 | 8.57e-07 | (-0.0000, -0.0000) | 0.176 | 2 | (-1,-1) |
| -(0.8900,0.7803) | 0.281 | 1.55e-10 | (-0.0000, -0.0000) | 0.175 | 18.3523 | -(0.8082,1.7803) |
| (-0.4612,0.2451) | 0.297 | 1.02e-12 | (-0.0000, -0.0000) | 0.181 | 16.9875 | -(0.1963,0.7549) |
| (0.2137,-0.0280) | 0.275 | 1.45e-11 | (0.0000,0.0000) | 0.173 | 21.3182 | (0.4754,-1.0280) |
| (0.7826,0.5242) | 0.234 | 2.73e-12 | (0.0000,0.0000) | 0.194 | 23.8254 | (0.8745,-0.4758) |
| -(0.0871,0.9630) | 0.343 | 1.49e-17 | (-0.0000,-0.0000) | 0.185 | 11.6507 | (0.2095,-1.9630) |
| (0.6428,-0.1106) | 0.203 | 1.29e-11 | (0.0000, 0.0000) | 0.189 | 11.8808 | (0.7868,-1.1106) |
| (0.2309,0.5839) | 0.312 | 5.71e-14 | (0.0000,- 0.0000) | 0.193 | 39.0324 | (0.4893,-0.4161) |
| (0.8436,0.4764) | 0.219 | 2.73e-12 | (0.0000,-0.0000) | 0.127 | 22.5167 | (0.9110,-0.5236) |
| -(0.6475,0.1886) | 0.203 | 2.68e-11 | (-0.0000, 0.0000) | 0.179 | 27.1169 | -(0.4390,1.1886) |
| (0.8709,0.8338) | 0.235 | 6.36e-11 | (0.0000,0.0000) | 0.182 | 6.8946 | (0.9270,-0.1662) |

The numerical simulation showed that any equilibrium points of the neural network of augmented smoothing were solutions for the problem 1 and 2, and our neural network can find globally and rapidly solutions of the nonlinear $L_1$-norm problem. The popular neural networks based on gradient system fail to find any solution of the nonconvex problems.

## 5. Concluding Remarks

Our research in this paper focuses on the augmented smoothing and the gradient-based neural network for the nonlinear $L_1$-norm problem (1). First, we introduce a new smoothing function of $L_1$-norm with an augmented regularization and study some properties of the smoothing function. Second, we consider stability properties of a gradient-based neural network model to minimize the smoothing function.

Our assumptions are necessary, because the Lipschitz continuity of gradient of the smoothing function guarantees the existence and uniqueness of the solution of the neural network.

We illustrated the advantage of our smoothing neural network model through numerical simulations for two difficult nonconvex problems and believe that similar conclusions would hold for different augmented smoothing gradient-based neural networks for nonlinear $L_1$-norm minimization problems.

## References

1. Y. S. Xia, C. Y. Sun, and W. X. Zheng, Discrete-time neural network for fast solving large linear $L_1$ estimation problems and its application to image restoration, *IEEE Trans. Neural Networks & Learning Systems*, 23(5), 812-820, 2012
2. Y. S. Xia, A compact cooperative recurrent neural network for computing general constrained $L_1$ norm estimators," *IEEE Trans. Signal Process.*, 57(9), 3693–3697, 2009
3. L. V. Ferreira et al., Solving systems of linear equations via gradient systems with discontinuous right-hand sides: Application to LS-SVM, *IEEE Trans. Neural Networks*, 16 (2), 501-505, 2005
4. L. Gao, Numerical algorithms for nonlinear $L_p$-norm problem and its extreme case. *J. Comput. Appl. Math.*, 129, 139-150, 2001
5. J. Sun , Solving the discrete $L_p$-approximation problem by a method of centers. *J. Comput. Appl. Math.*, 129, 63-76, 2001
6. Q. Han et al., Stability Analysis of Gradient-Based Neural Networks for Optimization Problems, *J. Global Optim.*, 19, 363-381, 2001
7. Z. Heszberger, J. Biro, An optimization neural network model with time-dependent and lossy dynamics, *Neurocomputing*, 48, 53-62, 2002
8. J. M. Peng, Z.H. Lin, A non-interior continuation method for generalized linear complementarity problem, *Math. Program.*, 86, 533-563, 1999
9. B. Yu et al., The aggregate constraint homotopy method for nonconvex nonlinear programming. *Nonlinear Analysis*, 45, 839- 847, 2001
10. M. G. Sklar, R. D. Armstrong, Lagrange method for large-scale minimal absolute estimation, *Computer & OR*. 20, 83-93, 1993
11. N. N. Abdelmalek, Efficient solution algorithm of discrete linear $L_1$-approximation problem, *Math.Comput.*,29, 844-850, 1975
12. I. Barrodale, F.D.K.Roberts, $L_1$-linear approximation and improved algorithm, *SIAM J.Num.Anal.*,10, 839-848, 1973